\documentclass[11pt]{article}

\usepackage{amsfonts,amssymb,amsmath,latexsym,indentfirst}

\usepackage[notcite,notref]{showkeys}

\newtheorem{prop}{Proposition}[section]
\newtheorem{rema}{Remark}[section]

\newtheorem{lemm}{Lemma}[section]
\newtheorem{theo}{Theorem}[section]

\newcommand{\R}{\ensuremath{{\mathbb{R}} }}

\newcommand{\peq}{\hspace*{0.10in}}
\newcommand{\peqq}{\hspace*{0.05in}}

\newcommand{\quebra}{\hspace*{-0.25in}}

\newcommand{\fim}{\rightline{$\blacksquare$}}

\author{Luiz Gustavo  Farah
\\
Department of Mathematics\\
University of California, Santa Barbara\\
California, 93106\footnote{Permanent address: Department of Mathematics - ICEx/UFMG, CEP 31270-901, Belo Horizonte, MG, Brazil}}
\title{Global rough solutions to the critical generalized KdV equation
\footnotetext{The author was supported by CNPq-Brazil.}}
\date{}

\begin{document}
\maketitle

\begin{abstract}
We prove that the initial value problem (IVP) for the critical generalized KdV
equation $u_{t}+u_{xxx}+(u^5)_{x}=0$ on the real line is globally well-posed in
$H^{s}(\R)$ in $s>3/5$ with the appropriate smallness assumption on the initial data.
\end{abstract}

\section{Introduction}
In this work we consider the initial value problem (IVP) for the critical generalized KdV
equation
\begin{eqnarray}\label{KdV}
\left\{
\begin{array}{l}
u_{t}+u_{xxx}+(u^5)_{x}=0, \peq  x\in \R,\, t>0,\\
u(x,0)=u_0(x).
\end{array} \right.
\end{eqnarray}

From the point of view of physics this kind of problem appears, for example, in the study of waves on shallow water (see Korteweg and de Vries \cite{KdV95}).

Well-posednes for the Cauchy problem (\ref{KdV}) has been studied by many authors. We refer the reader to Kato \cite{K83} (and references therein) for the $H^s$ theory ($s > 3/2$) and Kenig, Ponce and Vega \cite{KPV6} for the $L^2$ theory. We should notice that the latter result is optimal in view of the work Birnir, Kenig, Ponce, Svanstedt, Vega \cite{BKPSV96}. The results in \cite{KPV6} also imply that solutions corresponding to small data $u_0\in L^2(\R)$, say
\begin{equation}\label{EPS}
\|u_0\|_{L^2}\leq \epsilon_0
\end{equation}
are global in time. Note that this global $L^2$ result is valid for real or complex solutions and for both signs of the nonlinearity (focusing or defocusing). This is due to the homogeneity of the equation (scaling argument) and not to the $L^2$ conserved quantity.

It is know that real solutions for the equation (\ref{KdV}) satisfy the following conserved quantities
\begin{equation}\label{MC}
Mass\equiv M(t)=\|u(t)\|_{L^2};
\end{equation}
and
\begin{equation}\label{EC}
Energy\equiv E(t)=\frac{1}{2}\|u_x(t)\|^2_{L^2}-\frac{1}{6}\|u(t)\|^6_{L^6}.
\end{equation}

On the other hand, Weinstein \cite{W} showed the following sharp Gagliardo-Nirenberg inequality for $v\in H^1(\R)$ and $Q(x)=[3c\hspace*{0.025in} \textrm{sech}^2(2\sqrt{c}x)]^{1/4}$ (the solitary wave solution of (\ref{KdV}))
\begin{equation}\label{GNI}
\|v\|^6_{L^6}\leq3\left(\dfrac{\|v\|_{L^2}}{\|Q\|_{L^2}}\right)^4\|v_x\|^2_{L^2}.
\end{equation}

This estimate combined with the conserved quantities (\ref{MC}) and (\ref{EC}) force the energy to be positive and gives an \textit{a priori} estimate in $H^1(\R)$ provided 
\begin{equation}\label{SA}
\|u_0\|_{L^2}<\|Q\|_{L^2}.
\end{equation}

The local theory in $H^1(\R)$ together with the \textit{a priori} estimate immediately yield global-in-time well-posedness of (\ref{KdV}) from data $u_0 \in H^1(\R)$ under the smallness assumption (\ref{SA}). This result was improved by Fonseca, Linares and Ponce \cite{FLP}, who proved global well-posedness in $H^s(\R)$ for $s>3/4$, assuming (\ref{SA}). The method of proof is based on the idea of Bourgain \cite{B1}-\cite{B5} of estimating separately the evolution of low frequencies and of high frequencies. Indeed, it is expected that $\epsilon_0$ in \eqref{EPS} to be equal to the size of the solitary wave solution of (\ref{KdV}) (i.e., $\epsilon_0=\|Q\|_{L^2}$) and this is an interesting open problem (see Linares and Ponce \cite{LP09} page 185).

Several interesting results have been obtained for solutions of IVP (\ref{KdV}). Merle \cite{M01} (see also Martel and Merle \cite{MM02a}) proved the existence of real-valued solutions of (\ref{KdV}) in $H^1(\R)$ corresponding to data $u_0\in H^1(\R)$ satisfying $\|u_0\|_{L^2}>\|Q\|_{L^2}$ that blow-up. There are also various results concerning instability of solitary wave solutions as well as the structure of the blow-up formation obtained by Martel and Merle \cite{MM01} and \cite{MM02}.

Our principal aim is to loosen the regularity requirements on the initial data which ensure global-in-time solutions for the IVP (\ref{KdV}). In this paper, we prove the following resut. 
\begin{theo}\label{T1}
The initial value problem (\ref{KdV}) is globally well-posed in $H^s(\R)$ for all \peqq $s>3/5$, assuming the smallness condition (\ref{SA}). Moreover the solution satisfies
\begin{equation}\label{pb}
\sup_{t\in[0,T]}\left\{\|u(t)\|^2_{H^{s}}\right\}\leq C(1+T)^{\frac{1-s}{5s-3}+}
\end{equation}
where the constant $C$ depends only on $s$ and $\|u_0\|_{H^{s}}$.
\end{theo}

Here we use the approach introduced by Colliander,Keel, Staffilani, Takaoka and Tao in \cite{CKSTT4}, called the $I$-method. We also explain why the refined approach introduced by the same authors in \cite{CKSTT6}, \cite{CKSTT3} and \cite{CKSTT5} can not be use to improve our global result stated in Theorem \ref{T1} (see Proposition 3.1 and Remarks 3.1-3.2 below).

Note that when $u_0 \in H^{s}(\R)$ with
$s<1$ in (\ref{KdV}), the energy \eqref{EC} could be infinite, and so the
conservation law (\ref{EC}) is meaningless. To overcome this
difficulty, we follow the $I$-method scheme and introduce a modified
energy functional which is also defined for less regular functions.
Unfortunately, this new functional is not strictly conserved, but we
can show that it is \textit{almost} conserved in time. When one is
able to control its growth in time explicitly this allows to iterate
a modified local existence theorem to continue the solution to any
time $T$. 


The plan of this paper is as follows. In the next section we introduce some notation and preliminaries. Section 3 describe the multilinear correction technique which generates modified energies. In Section 4 we prove the almost conserved law.  Section 5 contains the variant local well-posedness result and the proof of the global result stated in Theorem \ref{T1}.

\section{Notations and preliminaries}\label{S3}

We use $c$ to denote various constant depending on $s$. Given any positive numbers $a$ and $b$, the notation $a \lesssim b$ means that there exists a positive
constant $c$ such that $a \leq cb$. Also, we denote $a \sim b$ when, $a \lesssim b$ and
$b \lesssim a$. We use $a+$ and $a-$ to denote $a+\varepsilon$ and $a-\varepsilon$, respectively, for
arbitrarily small exponents $\varepsilon>0$.

We use $\|f\|_{L^p}$ to denote the $L^p(\R)$ norm and $L^q_tL^r_x$ to denote the mixed norm
\begin{equation*}
\|f\|_{L^q_tL^r_x}\equiv \left(\int \|f\|_{L^r}^q \right)^{1/q}
\end{equation*}
with the usual modifications when $q =\infty$.\\

We define the spatial Fourier transform of $f(x)$ by
\begin{equation*}
\hat{f}(\xi)\equiv\int_{\R}e^{-ix\xi}f(x)dx
\end{equation*}
and the spacetime Fourier transform $u(t, x)$ by
\begin{equation*}
\widetilde{u}(\tau,\xi)\equiv\int_{\R}\int_{\R}e^{-i(x\xi+t\tau)}u(t,x)dtdx.
\end{equation*}
Note that the derivative $\partial_x$ is conjugated to multiplication by $i\xi$ by the Fourier transform.

We shall also define $D$ and $J$ to be, respectively, the Fourier multiplieres with symbol $|\xi|$ and $\langle \xi \rangle = 1+|\xi|$. Thus, the Sobolev norms $H^s(\R)$ is given by
\begin{equation*}
\|f\|_{H^s}\equiv \|J^sf\|_{L^2_x}=\|\langle \xi \rangle^s\hat{f}\|_{L^2_{\xi}},
\end{equation*}
where $\langle \xi \rangle = 1+|\xi|$.

We also define the $X_{s,b}(\R\times\R)$ spaces on $\R\times\R$ by 
\begin{equation*}
\|u\|_{X_{s,b}(\R\times\R)}=\|\langle\tau+|\xi|^2\rangle^b\langle\xi\rangle^s \widetilde{F}\|_{L^{2}_{\xi,\tau}}.
\end{equation*}

These spaces were used to systematically study nonlinear dispersive wave problems
by Bourgain \cite{B}. Klainerman and Machedon \cite{KM} used similar ideas in their
study of the nonlinear wave equation. The spaces appeared earlier in the study of propagation of singularity in semilinear wave equation in the works \cite{RR}, \cite{Be83} of Rauch, Reed, and M. Beals.

For any interval $I$ we define the localized $X_{s,b}(I \times\R)$ spaces by 
\begin{equation*}
\|u\|_{X_{s,b}(I \times\R)}=\inf\left\{\|w\|_{X_{s,b}(\R \times\R)}:w(t)=u(t) \textrm{ on }  I\right\}.
\end{equation*}

We often abbreviate $\|u\|_{X_{s,b}}$ and $\|u\|_{X^{I}_{s,b}}$, respectively, for $\|u\|_{X_{s,b}(\R\times\R)}$ and $\|u\|_{X_{s,b}(I \times\R)}$.

We shall take advantage of the Strichartz estimate (see Kenig, Ponce and Vega \cite{KPV4}) 
\begin{equation}\label{L8}
\|u\|_{L^8_{x,t}}\lesssim \|u\|_{X_{0,\frac{1}{2}+}}
\end{equation}
which interpolate with the trivial estimate
\begin{equation*}
\|u\|_{L^2_{x,t}}\lesssim \|u\|_{X_{0,0}}
\end{equation*}
to give
\begin{equation}\label{L6}
\|u\|_{L^6_{x,t}}\lesssim \|u\|_{X_{0,\frac{1}{2}+}}.
\end{equation}

We also use 
\begin{equation*}
\|u\|_{L^{\infty}_tL^2_{x}}\lesssim \|u\|_{X_{0,\frac{1}{2}+}},
\end{equation*}
which together with Sobolev embedding gives
\begin{equation}\label{LI}
\|u\|_{L^{\infty}_{x,t}}\lesssim \|u\|_{X_{\frac{1}{2}+,\frac{1}{2}+}}.
\end{equation}

Interpolation between \eqref{LI} and \eqref{L8} give us
\begin{equation}\label{LP}
\|u\|_{L^p_{x,t}}\lesssim \|u\|_{X_{\alpha(p),\frac{1}{2}+}},
\end{equation}
where $p>8$ and $\alpha(p)=\left(\dfrac{1}{2}+\right)\left(\dfrac{p-8}{p}\right)$. 

We also have the following refined Strichartz estimate in the case of differing frequencies (cf. Bourgain \cite{B1} and Gr\"unrock \cite{Gr05}).

\begin{lemm}\label{L3}
Let $\psi_1, \psi_2 \in X_{0,\frac{1}{2}+}$ be supported on spatial frequencies $|\xi_i|\sim N_i$, $i=1,2$.
If $\max\{|\xi_1|,|\xi_2|\} \lesssim \min\left\{|\xi_1-\xi_2|,|\xi_1+\xi_2|\right\}$ for all
$\xi_i\in \textrm{supp}(\widehat{\psi}_i)$, $i=1,2$, then
\begin{equation}\label{GRU}
\|\psi_1D_x\psi_2\|_{L^2_{x,t}}\lesssim \|\psi_1\|_{X_{0,\frac{1}{2}+}}
 \|\psi_2\|_{X_{0,\frac{1}{2}+}}.
\end{equation}
\end{lemm}

\noindent\textbf{Proof of Lemma \ref{L3}:} This is an improved Strichartz estimate of the type considered in Bourgain \cite{B1} (see also Ozawa and Tsutsumi \cite{OT}). In fact, the desired estimate is contained in Lemma 1 of Gr\"unrock \cite{Gr05}. We present the short proof for the sake of completeness. It is enough to show that
\begin{equation*}
\|(e^{-it\partial_x^3} u)(e^{-it\partial_x^3}D_xv)\|_{L^2_{x,t}}\lesssim \|u\|_{L^2} \|v\|_{L^2}
\end{equation*}
for functions $\widehat{u}(\xi_1)$ and $\widehat{v}(\xi_2)$ with support in $|\xi_i|\sim N_i$, $i=1,2$.

We have
\begin{eqnarray*}
\|(e^{-it\partial_x^3} u)(e^{-it\partial_x^3}D_xv)\|_{L^2_{x,t}}
\end{eqnarray*}
\vspace*{-0.25in}
\begin{eqnarray*}
&=&\int\!\!\!\int\left(\int\!\!\!\int_{\ast} e^{it(\xi_1^3+\xi_2^3-\eta_1^3-\eta_2^3)} |\xi_2||\eta_2|\widehat{u}(\xi_1)\widehat{v}(\xi_2) \overline{\widehat{u}(\eta_1)} \hspace*{0.05in}\overline{\widehat{v}(\eta_2)}d\xi_1d\eta_1\right)d\xi dt\\
&=&\int\left(\int\!\!\!\int_{\ast} \delta(P(\eta_1))|\xi_2||\eta_2|\widehat{u}(\xi_1)\widehat{v}(\xi_2) \overline{\widehat{u}(\eta_1)} \hspace*{0.05in}\overline{\widehat{v}(\eta_2)}d\xi_1d\eta_1\right)d\xi
\end{eqnarray*}
where $\ast$ denotes integration over $\xi=\xi_1+\xi_2=\eta_1+\eta_2$ and
\begin{equation*}
P(\eta_1)=\eta_1^3+\eta_2^3-\xi_1^3-\xi_2^3-=3\xi(\eta_1^2-\xi_1^2+\xi(\xi_1-\eta_1)).
\end{equation*}

Note that $P(\eta_1)$ has roots $\eta_1=\xi_1$ and $\eta_1=\xi-\xi_1$. Now, using the
well-known identity $\delta(g(x))=\sum_n\frac{\delta(x-x_n)}{|g^{'}(x_n)|}$, where the sum is
taken over all simple zeros of $g$, we obtain
\begin{eqnarray*}
\|(e^{-it\partial_x^3} u)(e^{-it\partial_x^3}D_xv)\|_{L^2_{x,t}}
\end{eqnarray*}
\vspace*{-0.25in}
\begin{eqnarray*}
&\lesssim&\int\left(\int \dfrac{|\xi_2|^2\widehat{u}(\xi_1) \overline{\widehat{u}(\xi_1)}\widehat{v}(\xi-\xi_1) \overline{\widehat{v}(\xi-\xi_1)}}{|\xi_1+\xi_2||\xi_1-\xi_2|} d\xi_1\right)d\xi\\
&&+\int\left(\int \dfrac{|\xi_1|^2\widehat{u}(\xi_1) \overline{\widehat{u}(\xi-\xi_1)}\widehat{v}(\xi-\xi_1) \overline{\widehat{v}(\xi_1)}}{|\xi_1+\xi_2||\xi_1-\xi_2|} d\xi_1\right)d\xi\\
&\lesssim&\|u\|_{L^2} \|v\|_{L^2},
\end{eqnarray*}
where in last inequality we have used the fact that 
$$\max\{|\xi_1|,|\xi_2|\} \lesssim \min\left\{|\xi_1-\xi_2|,|\xi_1+\xi_2|\right\}.$$\\
\fim

\begin{rema}
Note that the relation $|\xi_2|\gg|\xi_1|$ implies the hypothesis of the above lemma. This exactly the frequency assumption made in  Bourgain \cite{B1} for the Schr\"odinger equation.
\end{rema}

We now give some useful notation for multilinear expressions. If $n \geq 2$ is an even integer, we define a (spatial) $n$-multiplier to be any function $M_n(\xi_1,\dots,\xi_n)$ on the hyperplane
$$\Gamma_n\equiv\{ (\xi_1,\dots,\xi_n)\in \R^n : \xi_1+\cdots+\xi_n=0\},$$
which we endow with the standard measure $\delta(\xi_1+\cdots+\xi_n)$, where $\delta$ is the Dirac delta.

If $M_n$ is an $n$-multiplier and $f_1, \dots , f_n$ are functions on $\R$, we define the n-linear functional $\Lambda_n(M_n; f_1, \dots , f_n)$ by
$$\Lambda_n(M_n; f_1, \dots , f_n)\equiv \int_{\Gamma_n}M_n(\xi_1,\dots,\xi_n)\prod_{j=1}^n\widehat{f_j}(\xi_j).$$

We will often apply $\Lambda_n$ to $n$ copies of the same function $u$ in which case the dependence upon $u$ may be suppressed in the notation: $\Lambda_n(M_n; u, \dots , u)$ may simply be written $\Lambda_n(M_n)$.

If $M_n$ is symmetric, then is a symmetric k-linear functional $\Lambda_n(M_n)$.

As an example, suppose that $u$ is an $\R$-valued function. By Plancherel, we can rewrite the energy (\ref{EC}) in terms of $n$-linear functionals as
\begin{eqnarray*}
E(t)=-\dfrac{1}{2}\Lambda_2(\xi_1\xi_2)-\dfrac{1}{6}\Lambda_6(1).
\end{eqnarray*}

The time derivative of a symmetric n-linear functional can be calculated explicitly if we assume that the function $u$ satisfies a particular PDE. The following statement may be directly verified by using the critical generalized KdV equation (\ref{KdV}).
\begin{prop}
Suppose u satisfies the the critical generalized KdV equation (\ref{KdV}) and that $M_n$ is a symmetric $n$-multiplier. Then  
\begin{equation}\label{LD}
\dfrac{d}{dt}\Lambda_n(M_n)=\Lambda_n(M_n\alpha_n)-in\Lambda_{n+4}(M_n(\xi_1,\dots,\xi_{n-1},\xi_n+\cdots+\xi_{n+4})(\xi_n+\cdots+\xi_{n+4})),
\end{equation}
where $\alpha_n\equiv i(\xi_1^3+\cdots+\xi_n^3)$.
\end{prop}

\section{Modified energy functional}

As we mention in the introduction, we follow the ``almost conservation law" scheme introduced in Colliander, Keel, Staffilani, Takaoka and Tao \cite{CKSTT6}-\cite{CKSTT2}. To this end, we introduced a substitute notion of ``energy" that could be defined for less regular functions and that had very low increment in time.
Given $s<1$ and a parameter $N\gg 1$, define the multiplier operator $I_N:H^s\rightarrow H^1$ such that
\begin{equation*}
\widehat{I_Nf(\xi)}\equiv m_N(\xi)\widehat{f}(\xi),
\end{equation*}
where the multiplier $m_N(\xi)$ is smooth, radially symmetric, nondecreasing in $|\xi|$ and
\begin{eqnarray*}
m_N(\xi)=\left\{
\begin{array}{l l }
1&, \textrm{ if } |\xi|\leq N,\\
\left(\dfrac{N}{|\xi|}\right)^{1-s}&, \textrm{ if } |\xi|\geq 2N.
\end{array} \right.
\end{eqnarray*}

To simplify the notation, we omit the dependence of $N$ in $I_N$ and denote it only by $I$. Note that the operator $I$ is smoothing of order $1-s$. Indeed we have
\begin{equation}\label{smo}
\|u\|_{X_{s_0,b_0}}\leq c\|Iu\|_{X_{s_{0}+1-s,b_0}}\leq cN^{1-s}\|u\|_{X_{s_0,b_0}},
\end{equation}
for any $s_0,b_0\in \R$.

Our substitute energy will be defined by $E^1(u)=E(Iu)$. Obviously this energy makes sense even if $u$ is only in $H^s(\R)$. Thus,  in terms of $n$-linear functionals we have
\begin{equation}\label{ME1}
E^1(u)=-\dfrac{1}{2}\Lambda_2(m_1\xi_1m_2\xi_2)-\dfrac{1}{6}\Lambda_6(m_1\dots m_6),
\end{equation}
where $m_j=m(\xi_j)$.

We can think about $E^1(u)$ as the first generation of a family of modified energies. We also define the second energy
\begin{equation}\label{ME2}
E^2(u)=-\dfrac{1}{2}\Lambda_2(m_1\xi_1m_2\xi_2)-\dfrac{1}{6}\Lambda_6(M_6(\xi_1,\dots, \xi_6)).
\end{equation}

Thus, using the derivation law (\ref{LD}), we obtain
\begin{eqnarray*}
\frac{d}{dt}E^2(u)&=&-\dfrac{1}{2}\left[\Lambda_2(m_1\xi_1m_2\xi_2\alpha_2)-2i\Lambda_6(m_1\xi_1m(\xi_2+\cdots+\xi_6)(\xi_2+\cdots+\xi_{6})^2)\right]\\
&&-\dfrac{1}{6}\left[\Lambda_6(M_6\alpha_6)-6i\Lambda_{10}(M_6(\xi_1,\dots,\xi_{5},\xi_6+\cdots+\xi_{10})(\xi_6+\cdots+\xi_{10}))\right]\\
&=&-\dfrac{1}{2}\Lambda_2(m_1\xi_1m_2\xi_2\alpha_2)+\\
&&+\dfrac{i}{6}\Lambda_6(M_6(\xi_1^3+\dots+\xi_6^3)-(m_1^2\xi_1^3+\dots+m_6^2\xi_6^3))+\\
&&+i\Lambda_{10}(M_6(\xi_1,\dots,\xi_{5},\xi_6+\cdots+\xi_{10})(\xi_6+\cdots+\xi_{10}))
\end{eqnarray*}
where in the last equality we have used the identity $\xi_1+\cdots+\xi_6=0$ and symmetrizing.

Note that picking 
\begin{equation*}
M_6(\xi_1,\dots,\xi_{6})=\dfrac{m_1^2\xi_1^3+\dots+m_6^2\xi_6^3}{\xi_1^3+\dots+\xi_6^3}
\end{equation*}
we can force $\Lambda_6$ to be zero. Unfortunately the multiplier $M_6$ is not well defined in the set $\Gamma_6$. In fact, given $N\gg 1$, we can find numbers $\xi_1,\dots,\xi_6$ such that the denominator of $M_6$ is zero and the numerator is different from zero. This is the content of the following proposition.

\begin{prop}\label{P1}
Let $1/2<s<1$, there exist numbers $\xi_1,\dots,\xi_6$ such that
\begin{itemize}
\item[$(i)$] 
$\left\{ 
\begin{array}{l }
\xi_1+\cdots+\xi_6= 0;\\
\xi_1^3+\cdots+\xi_6^3=0.
\end{array} \right.$ 
\item[$(ii)$] $m_1^2\xi_1^3+\cdots+m_6^2\xi_1^3\neq 0.$
\end{itemize}
\end{prop}

\noindent\textbf{Proof. } There are several ways to find such numbers. Here, we only left to the reader the verification that the numbers $\xi_1=\xi_2=-k$, $\xi_3=-8k$, $\xi_4=\left(5+\dfrac{2\sqrt{55}}{5}\right)k$, $\xi_5=\left(5-\dfrac{2\sqrt{55}}{5}\right)k$ and $\xi_6=0$, where $k \gg N$, satisfy the relations $(i)$ and $(ii)$.\\
\fim

\begin{rema}
A similar conclusion can be made for the 3g-KdV equation
\begin{eqnarray*}
\left\{
\begin{array}{l}
u_{t}+u_{xxx}+(u^4)_{x}=0, \peq  x\in \R,\, t>0,\\
u(x,0)=u_0(x).
\end{array} \right.
\end{eqnarray*}

One can show that the refined approach proposed in \cite{CKSTT6}, \cite{CKSTT3} and \cite{CKSTT5} also does not work for this equation. In fact, as far as we know, the best rough global result up to now is given in Gr\"unrock, Panthee and Silva \cite{GPS07}, where the authors used the $I$-method in its first version.
\end{rema}

\begin{rema}
In general, we have the following. Define
\begin{equation}\label{Mj}
M_j(\xi_1,\dots,\xi_{j})=\dfrac{m_1^2\xi_1^3+\dots+m_j^2\xi_j^3}{\xi_j^3+\dots+\xi_j^3}.
\end{equation}

When $j=3,4$ the arithmetic facts (see, for example, Fefferman \cite{F73} and Bourgain \cite{B})
\begin{equation}\label{AF1}
\xi_1+\xi_2+\xi_3=0\Longrightarrow \xi_1^3+\xi_2^3+\xi_3^3=3\xi_1\xi_2\xi_3
\end{equation}
and
\begin{equation}\label{AF2}
\xi_1+\xi_2+\xi_3+\xi_4=0\Longrightarrow \xi_1^3+\xi_2^3+\xi_3^3+\xi_4^3=3(\xi_1+\xi_2)(\xi_1+\xi_3)(\xi_1+\xi_4)
\end{equation}
imply that the numerator must be zero if the denominator vanish in \eqref{Mj}. This fact was observed by  Colliander, Keel, Staffilani, Takaoka and Tao \cite{CKSTT3}, where the authors used the $I$-method in his refined approach to obtain sharp global well-posedness results for the KdV and modified KdV.

However, when $j\geq5$, due to the lack of an identity similar to \eqref{AF1}-\eqref{AF2}, one can prove an analogous result as the one stated in Proposition \ref{P1}. This implies that the multiplier $M_j$ is not well defined in this case. 
 
\end{rema}

Therefore, throughout this paper we will work only with the first modified energy (\ref{ME1}). Again, using the derivation law (\ref{LD}) and symmetrizing we have
\begin{eqnarray*}
\frac{d}{dt}E^1(u)&=&-\dfrac{i}{2}\Lambda_2(m_1\xi_1m_2\xi_2(\xi_1^3+\xi_2^3))+\\
&&+\dfrac{i}{6}\Lambda_6(m_1\dots m_6(\xi_1^3+\dots+\xi_6^3)-(m_1^2\xi_1^3+\dots+m_6^2\xi_6^3))+\\
&&+i\Lambda_{10}(m_1\dots m_5m(\xi_6+\cdots+\xi_{10})(\xi_6+\cdots+\xi_{10}))
\end{eqnarray*}

Observe that if $m=1$, the $\Lambda_6$ term vanish trivially. On the other hand, the terms $\Lambda_2$ and $\Lambda_{10}$ is also zero, since we have the restriction $\xi_1+\xi_2=0$ in the first and symmetrization in the later. This reproduces the Fourier proof of the energy conservation (\ref{EC}).

As one particular instance of the above computations and the Fundamental Theorem of Calculus, we have\\
\begin{eqnarray}\label{EMC1}
E^1(u)(t)-E^1(u)(0)=\int_0^t\frac{d}{dt}E^1(u)(t')dt'=
\end{eqnarray}
\vspace*{-0.25in}
\begin{eqnarray*}
&=&\dfrac{i}{6}\int_0^t\Lambda_6(m_1\dots m_6(\xi_1^3+\dots+\xi_6^3)-(m_1^2\xi_1^3+\dots+m_6^2\xi_6^3))(t')dt'\\
&&+i\int_0^t\Lambda_{10}(m_1\dots m_5m(\xi_6+\cdots+\xi_{10})(\xi_6+\cdots+\xi_{10}))(t')dt'.
\end{eqnarray*}

Most of our arguments here consist in  showing that the quantity $E^1(u)$ is \textit{almost} conserved
in time.

\section{Almost conservation law}

This section presents a detailed analysis of the expression (\ref{EMC1}). The analysis identifies some cancellations in the pointwise upper bound of some multipliers depending on the relative size of the frequencies envolved. Our aim is to prove the following almost conservation property.
\begin{prop}\label{p4.1}
Let $s>1/2$, $N\gg 1$ and $u\in H^s(\R)$ be a solution of \eqref{KdV} on $[T, T+\delta]$ such that $Iu\in H^1(\R)$. Then the following estimate holds
\begin{equation}\label{CC}
\left|E^1(u)(T+\delta)-E^1(u)(T)\right|\leq N^{-2+}\left( \left\|Iu\right\|_{X^{\delta}_{1,\frac{1}{2}+}}^6+\left\|Iu\right\|_{X^{\delta}_{1,\frac{1}{2}+}}^{10}\right).
\end{equation}
\end{prop}

\begin{rema}
The exponent $-2+$ on the right hand side of (\ref{CC}) is directly tied to the restriction $s > 3/5$ in our main theorem. If one could replace the increment $N^{-2+}$ by $N^{-\alpha+}$ for some $\alpha>0$ the argument we give in Section \ref{Global} implies global well-posedness of (\ref{KdV}) for all $s > 3/\alpha+3$. 
\end{rema}

\noindent\textbf{Proof. } We start with the estimative of the $\Lambda_6$ term. Instead of estimate each multilinear expression separately, we shall exploit some cancellation between the two multipliers. Using the fact that $\xi_1+\cdots+\xi_6=0$ this term can be rewritten as
\begin{equation*}
\Lambda_6(m_1\dots m_6(\xi_1^3+\dots+\xi_6^3)-(m_1^2\xi_1^3+\dots+m_6^2\xi_6^3))
\end{equation*}
\begin{equation*}
=-6\int_{\ast} \left( 1-\frac{m(\xi_2+\cdots+\xi_6)}{m(\xi_2)\cdots m(\xi_6)}\right) \xi_1^3\widehat{{Iu(\xi_1)}}
\cdots
 \widehat{Iu(\xi_6)},
\end{equation*}
where $\ast$ denotes integration over $\xi_1+\cdots+\xi_6=0$.

Therefore, our aim is to obtain the following inequality
\begin{equation*}
\mathbf{Term}\leq N^{-2+}\prod_{i=1}^6\left\|I\phi_i\right\|_{X^{\delta}_{1,\frac{1}{2}+}},
\end{equation*}
where
\begin{equation*}
\mathbf{Term}\equiv \left|\int_0^{\delta}\!\!\!\int_{\ast}  \left( 1-\frac{m(\xi_2+\cdots+\xi_6)}{m(\xi_2)\cdots m(\xi_6)}\right)\xi_1^3 \widehat{{I\phi_1(\xi_1)}}
 \cdots
 \widehat{I\phi_6(\xi_6)}\right|
\end{equation*}
and $\ast$ denotes integration over $\sum_{i=1}^6\xi_i=0$.

We estimate $\mathbf{Term}$ as follows. Without loss of generality, we assume the Fourier transforms of
all these functions to be nonnegative. First, we bound the symbol in the parentheses pointwise
in absolute value, according to the relative sizes of the frequencies involved. After that, the
remaining integrals are estimate using Plancherel formula, H\"older's inequality and Lemma
\ref{L3}. To sum over the dyadic pieces at the end we need to have extra factors
$N_j^{0-}$, $j=1,\dots,6$, everywhere.

We decompose the frequencies $\xi_j$, $j=1,\dots,6$, into dyadic blocks $N_j$. By the symmetry
of the multiplier
\begin{equation}\label{MULT}
1-\frac{m(\xi_2+\cdots+\xi_6)}{m(\xi_2)\cdots m(\xi_6)}
\end{equation}
in $\xi_2$, \dots, $\xi_6$, we may assume that
\begin{equation*}
N_2\geq \cdots \geq N_6.
\end{equation*}

Moreover, we can assume
$N_2 \gtrsim N$, because otherwise the symbol is zero. The condition $\sum_{i=1}^{6}\xi_i=0$
implies $N_1\lesssim N_2$. We split the different frequency interaction into several cases,
according to the size of the parameter $N$ in comparison to the $N_i$'s.\\

\quebra \textbf{Case $A$: }$N_2\gtrsim N\gg N_3\geq \cdots \geq N_6$.\\

The condition $\sum_{i=1}^{6}\xi_i=0$ implies $N_1\sim N_2$. By mean value theorem,
\begin{equation*}
\left|\frac{m(\xi_2)-m(\xi_2+\cdots+\xi_6)}{m(\xi_2)}\right|\lesssim \frac{\left|\nabla m(\xi_2)(\xi_3+\cdots+\xi_6)\right|}{m(\xi_2)}\lesssim\frac{N_3}{N_2}.
\end{equation*}

Therefore, Lemma \ref{L3} and (\ref{LI}) imply that
\begin{eqnarray*}
{\mathbf{Term}} \!\!\!&\lesssim& \!\!\! \frac{N_1^3N_3}{N_2}\left\|I\phi_1I\phi_3\right\|_{L^2(\R\times[0,\delta])} \left\|I\phi_2I\phi_4\right\|_{L^2(\R\times[0,\delta])} \left\|I\phi_5\right\|_{L^{\infty}}\left\|I\phi_6\right\|_{L^{\infty}}\\
&\lesssim&  \!\!\!\frac{N_1^3N_3}{N_2N_1N_2N_1N_2\langle N_3 \rangle\langle N_4 \rangle \langle N_5 \rangle^{1/2-} \langle N_6 \rangle^{1/2-}} \prod_{i=1}^{6} \|I\phi_i\|_{X^{\delta}_{1,\frac{1}{2}+}}\\
&\lesssim& \!\!\! N^{-2+}N_{max}^{0-}\prod_{i=1}^{6} \|I\phi_i\|_{X^{\delta}_{1,\frac{1}{2}+}}.
\end{eqnarray*}
\textbf{Case $B$: }$N_2\gg N_3\gtrsim N$ and $N_3\geq \cdots \geq N_6$.\\

In this case we also have $N_1\sim N_2$. We bound the multiplier (\ref{MULT}) by
\begin{equation}\label{MULT2}
\left|1-\frac{m(\xi_2+\cdots+\xi_6)}{m(\xi_2)\cdots m(\xi_6)}\right|\lesssim \frac{m(\xi_1)}{m(\xi_2) \cdots m(\xi_6)}.
\end{equation}

Therefore, since $m(N_1)\sim m(N_2)$, applying Lemma \ref{L3} and (\ref{LI}) we have
\begin{eqnarray*}
{\mathbf{Term}} \!\!\!&\lesssim&  \!\!\!\frac{N_1^3}{m(N_3)\cdots m(N_6)}\left\|I\phi_1I\phi_3\right\|_{L^2(\R\times[0,\delta])} \left\|I\phi_2I\phi_4\right\|_{L^2(\R\times[0,\delta])}  \left\|I\phi_5\right\|_{L^{\infty}}\left\|I\phi_6\right\|_{L^{\infty}}\\
&\lesssim& \!\!\! \frac{N_1^3} {m(N_3)\cdots m(N_6)N_1N_2N_1N_2N_3\langle N_4 \rangle \langle N_5 \rangle^{1/2-} \langle N_6 \rangle^{1/2-}}\prod_{i=1}^{6} \|I\phi_i\|_{X^{\delta}_{1,\frac{1}{2}+}}\\
&\lesssim& \!\!\! \frac{1} {m(N_3)N_3m(N_4)\langle N_4\rangle m(N_5)\langle N_5 \rangle^{1/2-} m(N_6)\langle N_6 \rangle^{1/2-}N_2}\prod_{i=1}^{6} \|I\phi_i\|_{X^{\delta}_{1,\frac{1}{2}+}}\\
&\lesssim& \!\!\! N^{-2+}N_{max}^{0-} \prod_{i=1}^{6} \|I\phi_i\|_{X^{\delta}_{1,\frac{1}{2}+}},
\end{eqnarray*}
where in the last inequality we use the fact that for any $p >0$ such that $p+s\geq 1$, the function
$m(x)x^p$ is increasing and $m(x)\langle x\rangle^p$ is bounded below, which implies
$m(N_3)N_3\gtrsim m(N)N=N$, $m(N_4)\langle N_4\rangle \gtrsim 1$ and $m(N_j)\langle N_j \rangle^{1/2-}\gtrsim 1$ for $j=5,6$.\\

\quebra \textbf{Case $C$: }$N_2\sim N_3\gtrsim N$ and $N_3\geq \cdots \geq N_4$.\\

The condition $\sum_{i=1}^{6}\xi_i=0$ implies $N_1\lesssim N_2$. We again bound the multiplier (\ref{MULT}) pointwise by (\ref{MULT2}). To obtain the decay $N^{-2+}$ we split this case into five subcases.\\

\textbf{Case $C.1$: } $N_4\gtrsim N$ and $N_4\ll N_3$.\\

From (\ref{MULT2}) and Lemmas \ref{L3}, we have that
\begin{equation*}
\begin{split}
{\mathbf{Term}} &\lesssim \frac{N_1^3m(N_1)}{m(N_2)\cdots m(N_6)} \left\|I\phi_2I\phi_4\right\|_{L^2(\R\times[0,\delta])} \left\|I\phi_3I\phi_5\right\|_{L^2(\R\times[0,\delta])}  \left\|I\phi_1\right\|_{L^{\infty}}\left\|I\phi_6\right\|_{L^{\infty}}\\
&\lesssim \frac{N_1^3m(N_1)} {m(N_2)\cdots m(N_6)N_2N_3N_2N_3N_4\langle N_5 \rangle \langle N_1 \rangle^{1/2-} \langle N_6 \rangle^{1/2-} } \prod_{i=1}^{6} \|I\phi_i\|_{X^{\delta}_{1,\frac{1}{2}+}}\end{split}
\end{equation*}
\begin{equation*}
\begin{split}
&\lesssim \! \frac{N_{max}^{0-}} {m(N_2)N_2^{3/4-}m(N_3)N_3^{3/4-}m(N_4)\langle N_4\rangle m(N_5)\langle N_5 \rangle m(N_6)\langle N_6 \rangle^{1/2-}}\prod_{i=1}^{6} \|I\phi_i\|_{X^{\delta}_{1,\frac{1}{2}+}}\\
&\lesssim \! N^{-2+}N_{max}^{0-} \prod_{i=1}^{6} \|I\phi_i\|_{X^{\delta}_{1,\frac{1}{2}+}}.
\end{split}
\end{equation*}

\textbf{Case $C.2$: } $N_4\gtrsim N$ and $N_3\sim N_4\gg N_5$.\\

Applying the same arguments as above
\begin{equation*}
\begin{split}
{\mathbf{Term}} &\lesssim \frac{N_1^3m(N_1)}{m(N_2)\cdots m(N_6)} \left\|I\phi_2I\phi_5\right\|_{L^2(\R\times[0,\delta])} \left\|I\phi_3I\phi_6\right\|_{L^2(\R\times[0,\delta])}  \left\|I\phi_1\right\|_{L^{\infty}}\left\|I\phi_4\right\|_{L^{\infty}}\\
&\lesssim \frac{N_1^3m(N_1)} {m(N_2)\cdots m(N_6)N_2N_3N_2N_3\langle N_5 \rangle \langle N_6 \rangle \langle N_1 \rangle^{1/2-} N_4^{1/2-} } \prod_{i=1}^{6} \|I\phi_i\|_{X^{\delta}_{1,\frac{1}{2}+}}\end{split}
\end{equation*}
\begin{equation*}
\begin{split}
&\lesssim \! \frac{N_{max}^{0-}} {m(N_2)N_2^{2/3-}m(N_3)N_3^{2/3-}m(N_4)N_4^{2/3-} m(N_5)\langle N_5 \rangle m(N_6)\langle N_6 \rangle^{1/2-}}\prod_{i=1}^{6} \|I\phi_i\|_{X^{\delta}_{1,\frac{1}{2}+}}\\
&\lesssim \! N^{-2+}N_{max}^{0-} \prod_{i=1}^{6} \|I\phi_i\|_{X^{\delta}_{1,\frac{1}{2}+}}.
\end{split}
\end{equation*}

\textbf{Case $C.3$:} $N_4\gtrsim N$ and $N_3\sim N_4\sim N_5$.\\

In view of \eqref{L6}, we have
\begin{equation*}
\begin{split}
{\mathbf{Term}} &\lesssim \frac{N_1^3m(N_1)}{m(N_2)\cdots m(N_6)} \prod_{i=1}^{6} \|I\phi_i\|_{L^6}\\
&\lesssim \frac{N_1^3m(N_1)} {m(N_2)\cdots m(N_6)\langle N_1 \rangle N_2N_3N_4N_5 \langle N_6 \rangle } \prod_{i=1}^{6} \|I\phi_i\|_{X^{\delta}_{1,\frac{1}{2}+}}\\
&\lesssim \frac{N_{max}^{0-}} {m(N_2)N_2^{1/2-}m(N_3)N_3^{1/2-}m(N_4)N_4^{1/2-} m(N_5)N_5^{1/2-} m(N_6)\langle N_6 \rangle}\prod_{i=1}^{6} \|I\phi_i\|_{X^{\delta}_{1,\frac{1}{2}+}}\\
&\lesssim N^{-2+}N_{max}^{0-} \prod_{i=1}^{6} \|I\phi_i\|_{X^{\delta}_{1,\frac{1}{2}+}}.
\end{split}
\end{equation*}

\textbf{Case $C.4$:} $N_4\ll N$ and $N_1\ll N_2$.\\

Again using the bound (\ref{MULT2}) and Lemma \ref{L3}, we have
\begin{equation*}
\begin{split}
{\mathbf{Term}} &\lesssim \frac{N_1^3m(N_1)}{m(N_2)\cdots m(N_6)}\left\|I\phi_2I\phi_1\right\|_{L^2(\R\times[0,\delta])} \left\|I\phi_3I\phi_4\right\|_{L^2(\R\times[0,\delta])}  \left\|I\phi_5\right\|_{L^{\infty}}\left\|I\phi_6\right\|_{L^{\infty}}\\
&\lesssim \frac{N_1^3m(N_1)} {m(N_2)\cdots m(N_6)N_2N_3\langle N_1 \rangle N_2N_3 \langle N_4 \rangle  \langle N_5 \rangle^{1/2-} \langle N_6 \rangle^{1/2-} } \prod_{i=1}^{6} \|I\phi_i\|_{X^{\delta}_{1,\frac{1}{2}+}}\\
&\lesssim \frac{N_{max}^{0-}} {m(N_2)N_2^{1-}m(N_3)N_3^{1-}m(N_4)\langle N_4 \rangle m(N_5)\langle N_5 \rangle^{1/2-} m(N_6)\langle N_6 \rangle^{1/2-}}\prod_{i=1}^{6} \|I\phi_i\|_{X^{\delta}_{1,\frac{1}{2}+}}\\
&\lesssim N^{-2+}N_{max}^{0-} \prod_{i=1}^{6} \|I\phi_i\|_{X^{\delta}_{1,\frac{1}{2}+}}.
\end{split}
\end{equation*}

\textbf{Case $C.5$: } $N_4\ll N$ and $N_1\sim N_2\sim N_3 \gtrsim N$.\\

In this case, we use an argument similar to the one used in Pecher \cite{P1} Proposition 5.1. Because
of $\sum_{i=1}^{6}\xi_i=0$, two of the large frequencies have different sign,
say, $\xi_1$ and $\xi_2$. Thus,
\begin{equation*}
|\xi_2|\leq |\xi_1-\xi_2|\leq 2|\xi_2|
\end{equation*}
and
\begin{equation*}
|\xi_1+\xi_2|= |\xi_3+\cdots+\xi_6|\sim|\xi_2|.
\end{equation*}

Therefore, using the bound \eqref{MULT2} and Lemma \ref{L3}, we obtain
\begin{equation*}
\begin{split}
{\mathbf{Term}} &\lesssim \frac{N_1^3m(N_1)}{m(N_2)\cdots m(N_6)}\left\|I\phi_2I\phi_1\right\|_{L^2(\R\times[0,\delta])} \left\|I\phi_3I\phi_4\right\|_{L^2(\R\times[0,\delta])}  \left\|I\phi_5\right\|_{L^{\infty}}\left\|I\phi_6\right\|_{L^{\infty}}\\
&\lesssim \frac{N_1^3m(N_1)} {m(N_2)\cdots m(N_6)N_2N_3N_1 N_2N_3 \langle N_4 \rangle  \langle N_5 \rangle^{1/2-} \langle N_6 \rangle^{1/2-} } \prod_{i=1}^{6} \|I\phi_i\|_{X^{\delta}_{1,\frac{1}{2}+}}\\
&\lesssim \frac{N_{max}^{0-}} {m(N_2)N_2^{1-}m(N_3)N_3^{1-}m(N_4)\langle N_4 \rangle m(N_5)\langle N_5 \rangle^{1/2-} m(N_6)\langle N_6 \rangle^{1/2-}}\prod_{i=1}^{6} \|I\phi_i\|_{X^{\delta}_{1,\frac{1}{2}+}}\\
&\lesssim N^{-2+}N_{max}^{0-} \prod_{i=1}^{6} \|I\phi_i\|_{X^{\delta}_{1,\frac{1}{2}+}}.
\end{split}
\end{equation*}

Now we turn to the estimate of the $\Lambda_10$ term. Before we start let us fix some notation. We write $N_1^{\ast}\geq N_2^{\ast}\geq N_3^{\ast}$ for the highest, second highest, third highest values of the frequencies $N_1,\cdots,N_6$. It is clear that 
\begin{equation}\label{B10}
|m_1\dots m_5m(\xi_6+\cdots+\xi_{10})(\xi_6+\cdots+\xi_{10})|\lesssim N_1^{\ast}.
\end{equation}

Again we perform a Littlewood-Paley decomposition of the ten functions $u$.\\

 \textbf{Case $A$: }$N_1^{\ast}\sim N_2^{\ast}\sim N_3^{\ast} \gtrsim N. $\\
 
 In view of \eqref{B10} and the fact that $m^3(N_1^{\ast})N_1^{\ast3-}\gtrsim N^3$, we have
 \begin{eqnarray*}
 \left|\int_T^{T+\delta}\Lambda_{10}(1)(t')dt'\right|&\lesssim& \dfrac{N^{\ast0-}_1}{N^{2-}}\int\!\!\!\int |JIu|^3|u|^7\\
 &\lesssim&\dfrac{N^{\ast0-}_1}{N^{2-}}\|JIu\|_{L^8}^3\|u\|^7_{56/5}\\
  &\lesssim&\dfrac{N^{\ast0-}_1}{N^{2-}}\|Iu\|_{X^{\delta}_{1,\frac{1}{2}+}}^3  \|u\|^7_{X^{\delta}_{\alpha(56/5),\frac{1}{2}+}},
 \end{eqnarray*}
where we have applied H\"older inequality, \eqref{L8} and \eqref{LP}. 

Note that $\alpha(56/5)=1/7+$. Therefore the inequality \eqref{smo} implies 
$$\|u\|_{X^{\delta}_{\alpha(56/5),\frac{1}{2}+}}\lesssim \|Iu\|_{X^{\delta}_{1,\frac{1}{2}+}},$$ 
for all $s>3/5$.

So, in this case
\begin{eqnarray*}
 \left|\int_T^{T+\delta}\Lambda_{10}(1)(t')dt'\right|&\lesssim& \dfrac{N^{\ast0-}_1}{N^{2-}}\|Iu\|_{X^{\delta}_{1,\frac{1}{2}+}}^{10}.
 \end{eqnarray*}
 
 \textbf{Case $B$: }$N_1^{\ast}\sim N_2^{\ast}\gg N_3^{\ast}$.\\
 
 Let $u_j\equiv u(N_j)$. Again, the inequality $m^2(N_1^{\ast})N_1^{\ast2-}\gtrsim N^2$ and \eqref{B10} implies that
 \begin{eqnarray*}
 \left|\int_T^{T+\delta}\Lambda_{10}(1)(t')dt'\right|&\lesssim& \dfrac{N^{\ast0-}_1}{N^{1-}}\|JIu_1u_3\|_{L^2}\|JIu_2\prod_{j=4}^{10}u_j\|_{L^2}\\
 &\lesssim& \dfrac{N^{\ast0-}_1}{N^{2-}}\|JIu_1\|_{L^2}\|u_3\|_{L^2}\|JIu_2\|_{L^8}\|u\|^7_{L^{56/3}}\\
  &\lesssim&\dfrac{N^{\ast0-}_1}{N^{2-}}\|Iu\|_{X^{\delta}_{1,\frac{1}{2}+}}^{10},
 \end{eqnarray*}
where we have applied H\"older inequality, \eqref{L8} and \eqref{LP} with $\alpha(56/3)=2/7+<3/5$.

This concludes the proof of Proposition \ref{p4.1}.\\
\fim

\section{Global theory}\label{Global}

Before proceed to the proof of Theorem \ref{T1} we need to establish a variant local well-posedness result for the following modified equation.
\begin{eqnarray}\label{MKdV}
\left\{
\begin{array}{l}
Iu_{t}+Iu_{xxx}+I(u^5)_{x}=0, \peq  x\in \R,\, t>0,\\
Iu(x,0)=Iu_0(x).
\end{array} \right.
\end{eqnarray}

Clearly if $Iu\in H^1(\R)$ is a solution of \eqref{MKdV}, then $u\in H^s(\R)$ is a solution of \eqref{KdV} in the same time interval. Therefore, we need to prove that, in fact, the above modified equation has a global solution.

The crucial nonlinear estimate for the local existence is given in the next lemma.
\begin{lemm}\label{NLEL}
For $s>1/2$ and $-\frac{1}{2}<b'\leq0<\frac{1}{2}<b$ we have
\begin{equation}\label{NLE}
\|\partial_x(\prod_{j=1}^5u_j)\|_{X_{s,b'}}\lesssim \prod_{j=1}^5\|u_j\|_{X_{s,b}}.
\end{equation}
\end{lemm}

\noindent\textbf{Proof. } For the sake of simplicity we will assume that $b'=0$ (clearly the general result is implied by this particular case). By definition
\begin{equation*}
\|\partial_x(\prod_{j=1}^5u_j)\|_{X_{s,0}}= c\left\|\xi\langle\xi\rangle^s \int_{\ast}\prod_{j=1}^{5} \widetilde{u_j}(\xi_j,\tau_j)\right\|^5_{L^2_{\xi,\tau}},
\end{equation*}
where $\ast$ denotes integration over $\sum(\xi_j,\tau_j)=(\xi,\tau).$

Again, we split the domain of integration according to the relative sizes of the spacial frequencies involved. By symmetry we may assume
\begin{equation*}
N_1\geq \cdots \geq N_6.
\end{equation*}

We will consider the following three regions.
\begin{eqnarray*}
A&=&\{ |\xi_1|\leq 1\};\\
B&=&\{ |\xi_1|\geq 1 \textrm{ and } |\xi_6|\leq |\xi_1|/2 \};\\
C&=&\{|\xi_1|\geq 1 \textrm{ and } |\xi_6|\geq |\xi_1|/2 \}.
\end{eqnarray*} 

In region $A$ we have $|\xi\langle\xi\rangle^{s}|\lesssim 1$. Therefore, by inequalities \eqref{L8} and \eqref{LI}, we conclude
\begin{eqnarray*}
\|\partial_x(\prod_{j=1}^5u_j)\|_{X_{s,0}}&\lesssim&  \prod_{j=1}^4\|u_j\|_{L_{x,t}^8} \|u_5\|_{L_{x,t}^{\infty}}\\
&\lesssim& \prod_{j=1}^5\|u_j\|_{X_{s,b}}.
\end{eqnarray*}

In region $B$ we have $|\xi_1| \lesssim \min\left\{|\xi_1-\xi_5|,|\xi_1+\xi_5|\right\}$. Applying Lemma \ref{L3} and inequalities \eqref{L8} and \eqref{LI}, we obtain
\begin{eqnarray*}
\|\partial_x(\prod_{j=1}^5u_j)\|_{X_{s,0}}&\lesssim& \|(D_xJ^su_1)u_5\|_{L_{x,t}^2} \prod_{j=2}^4\|u_j\|_{L_{x,t}^{\infty}}\\
&\lesssim& \prod_{j=1}^5\|u_j\|_{X_{s,b}}.
\end{eqnarray*}

In region $C$ we have $|\xi_1|\sim|\xi_j|$, for all $j=2,\dots,5$. Therefore $|\xi\langle\xi\rangle^{s}|\lesssim \langle\xi_1\rangle^{s}\prod_{j=2}^5\langle\xi_j\rangle^{1/4}$. Now, applying inequalities \eqref{L8} and \eqref{LP}, we have
\begin{eqnarray*}
\|\partial_x(\prod_{j=1}^5u_j)\|_{X_{s,0}}&\lesssim&  \|J^su_1\|_{L_{x,t}^8} \prod_{j=2}^5 \|J^{1/4}u_j\|_{L_{x,t}^{32/3}}\\
&\lesssim&  \|J^su_1\|_{X_{0,b}} \prod_{j=2}^5 \|J^{1/4}u_j\|_{X_{\alpha(32/3),b}}\\
&\lesssim& \prod_{j=1}^5\|u_j\|_{X_{s,b}},
\end{eqnarray*}
where we have used that $\alpha(32/3)+1/4=1/8+1/4<1/2$.\\
\fim
 
\begin{rema}
It should be interesting to prove inequality \eqref{NLE} for $s>0$ and the same assumptions on the parameters $b$ and $b'$. As a consequence, one can recover all the well know range of existence for the local theory in terms of the $X_{s,b}$ spaces.
\end{rema}

Applying the interpolation lemma (see \cite{CKSTT2}, Lemma 12.1) to \eqref{NLE} with $b'=0$ we obtain
\begin{equation*}
\|\partial_xI(\prod_{j=1}^5u_j)\|_{X_{1,0}}\lesssim \prod_{j=1}^5\|Iu_j\|_{X_{1,b}}.
\end{equation*}
where the implicit constant is independent of $N$. Now standard arguments invoking the contraction-mapping principle give the following variant of local well posedness.

\begin{theo}\label{t3.2}
Assume $s < 1$. Let $u_0 \in H^s(\R)$ be
given. Then there exists a positive number $\delta$ such
that the IVP (\ref{MKdV}) has a unique local solution $Iu\in C([0,\delta]:H^1(\R))$ such that
\begin{equation}\label{CONT}
\|Iu\|_{X^{\delta}_{1,\frac{1}{2}+}}\lesssim \|Iu_0\|_{H^1}.
\end{equation}

Moreover, the existence time can be estimated by
\begin{equation*}
\delta^{\frac{1}{2}-} \sim \dfrac{1}{\|Iu_0\|^4_{H^1}}.
\end{equation*}
\end{theo}

Now, we have all tools to proof our global result stated in Theorem \ref{T1}.\\

\noindent\textbf{Proof of Theorem \ref{T1}. }Let $u_0 \in H^s(\R)$ with $1/2< s <1$. Our goal is to construct a solution to \eqref{MKdV} (and therefore to \eqref{KdV}) on an arbitrary time interval $[0,T]$. We rescale the solution by writing $u_{\lambda}(x,t)=\lambda^{-1/2}u(x/\lambda,t/\lambda^3)$. We can easily check that $u(x,t)$ is a solution of \eqref{KdV} on the time interval $[0,T]$ if and only if $u_{\lambda}(x,t)$ is a solution to the same equation, with initial data $u_{0,\lambda}=\lambda^{-1/2}u_0(x/\lambda)$, on the time interval $[0,\lambda^3T]$. 

Since $|m(\xi)|\leq1$, in view of \eqref{MC}, a calculation shows that
\begin{equation*}
\|Iu_{\lambda}(t)\|_{L^2}\leq \|Iu_{0,\lambda}\|_{L^2}=\|Iu_0\|_{L^2}<\|Q\|_{L^2}.
\end{equation*}
where we have used the smallness condition \eqref{SA}.

Therefore, by the sharp Gagliardo-Niremberg \eqref{GNI}, we have
\begin{equation}\label{H1B}
\|\partial_xIu_{\lambda}(t)\|^2_{L^2}\lesssim E(Iu_{\lambda})(t).
\end{equation}
and
\begin{equation*}
E(Iu_{0,\lambda}) \lesssim \|\partial_xIu_{0,\lambda}\|^2_{L^2}.
\end{equation*}

On the other hand, using that $\dfrac{m(\xi)}{|\xi|^{s-1}}\lesssim N^{1-s}$, we obtain
\begin{equation}\label{IU0}
\|\partial_xIu_{0,\lambda}\|_{L^2}= \|m(\xi)|\xi|\widehat{u_{0,\lambda}}\|_{L^2}\leq N^{1-s}\||\xi|^s\widehat{u_{0,\lambda}}\|_{L^2}<\dfrac{N^{1-s}}{\lambda^s}\|u_0\|_{\dot{H}^s}.
\end{equation}

We use our variant local existence Theorem \ref{t3.2} on $[0, \delta]$, where $\delta^{\frac{1}{2}-}\sim \|Iu_{0,\lambda}\|^{-4}_{H^1}$ and conclude
\begin{equation}\label{IUV}
\|Iu_{\lambda}\|_{X^{\delta}_{1,\frac{1}{2}+}}\lesssim \|Iu_{0,\lambda}\|_{H^1}.
\end{equation}

The choice of the parameter $N=N(T)$ will be made later, but we select $\lambda$ now by requiring
\begin{eqnarray*}
\dfrac{N^{1-s}}{\lambda^s}\|u_0\|_{\dot{H}^s}<1\Longrightarrow \lambda \sim N^{\frac{1-s}{s}}.
\end{eqnarray*}

From now on, we drop the $\lambda$ subscript on $u$. By the almost conservation law stated in Proposition \ref{p4.1} and \eqref{IU0}-\eqref{IUV}, we have
\begin{equation*}
E^1(1)\leq E^1(0)+cN^{-2+}<1+cN^{-2+}<2.
\end{equation*}

We iterate this process $M$ times obtaining 
\begin{equation}\label{UB}
E^1(M)\leq E^1(0)+cMN^{-2+}<1+cMN^{-2+}<2,
\end{equation}
as long as $MN^{-2+}\lesssim 1$, which implies that the lifetime of the local results remains uniformly of size $1$. We take $M\sim N^{2-}$. This process extends the local solution to the time interval $[0,N^{2-}]$. Now, we choose $N=N(T)$ so that
$$N^{2-}>\lambda^3T\sim N^{3\left(\frac{1-s}{s}\right)}T\Longrightarrow N^{\frac{5s-3}{s}-}>T.$$
Therefore, if $s>\frac{3}{5}$ then $T$ can be taken arbitrarily large which conclude our global result.

Finally, we need to establish the polynomial bound (\ref{pb}). Undoing the scaling we have that
\begin{equation*}
\|\partial_xIu_{\lambda}(\lambda^3T_0)\|_{L^2}=\frac{1}{\lambda^2}\|\partial_xIu(T_0)\|_{L^2}.
\end{equation*}

Let $T_0\sim N^{\frac{5s-3}{s}-}$, therefore our uniform bound \eqref{UB} together with \eqref{smo}, \eqref{MC} and \eqref{H1B}  imply
\begin{eqnarray*}
\|u(T_0)\|^2_{H^{s}}\lesssim \|Iu(T_0)\|^2_{H^{1}}&\lesssim& \|Iu(T_0)\|^2_{L^{2}}+ \|\partial_xIu(T_0)\|^2_{L^{2}}\\
&\lesssim& \|u_0\|^2_{L^{2}}+ \lambda^2\|\partial_xIu_{\lambda}(\lambda^3T_0)\|^2_{L^{2}}\\
&\lesssim& \|u_0\|^2_{L^{2}}+ N^{2\left(\frac{1-s}{s}\right)}\\
&\lesssim& (1+T)^{\frac{1-s}{5s-3}+}
\end{eqnarray*}
\fim

\centerline{\textbf{Acknowledgment}}

This research was carried out when the author was visiting the Department of Mathematics of the University of California, Santa Barbara, whose hospitality is gratefully acknowledge.



E-mail: lgfarah@gmail.com.
\end{document}